\documentclass[12pt]{amsart}

\usepackage{amssymb,amsmath,xcolor,graphicx,xspace,colortbl,rotating} %
\usepackage[raggedrightboxes]{ragged2e} 
\usepackage{amsfonts}  
\usepackage{amsmath}  
\usepackage{amscd}  
\usepackage{amssymb}  
\usepackage{amssymb,amsmath,xcolor,graphicx,xspace,colortbl,rotating}  
\usepackage{graphicx}  
\usepackage{latexsym}  
\usepackage[parfill]{parskip}  
\usepackage[raggedrightboxes]{ragged2e}  
\usepackage[notcite,notref]{showkeys}  
\usepackage[all]{xy}  
\graphicspath{{TopK3.SWP_graphics/}{TopK3.SWP_tcache/}{TopK3.SWP_gcache/}}
\DeclareGraphicsExtensions{.pdf,.eps,.ps,.png,.jpg,.jpeg}
\graphicspath{{TopK2.SWP_graphics/}{TopK2.SWP_tcache/}{TopK2.SWP_gcache/}}
\DeclareGraphicsExtensions{.pdf,.eps,.ps,.png,.jpg,.jpeg}
\graphicspath{{TopK.SWP_graphics/}{TopK.SWP_tcache/}{TopK.SWP_gcache/}}
\DeclareGraphicsExtensions{.pdf,.eps,.ps,.png,.jpg,.jpeg}
\graphicspath{{Essai51.SWP_graphics/}{Essai51.SWP_tcache/}{Essai51.SWP_gcache/}}
\DeclareGraphicsExtensions{.pdf,.eps,.ps,.png,.jpg,.jpeg}
\graphicspath{{Essai49.SWP_graphics/}{Essai49.SWP_tcache/}{Essai49.SWP_gcache/}}
\DeclareGraphicsExtensions{.pdf,.eps,.ps,.png,.jpg,.jpeg}
\graphicspath{{Essai48.SWP_graphics/}{Essai48.SWP_tcache/}{Essai48.SWP_gcache/}}
\DeclareGraphicsExtensions{.pdf,.eps,.ps,.png,.jpg,.jpeg}
\graphicspath{{Essai47.SWP_graphics/}{Essai47.SWP_tcache/}{Essai47.SWP_gcache/}}
\DeclareGraphicsExtensions{.pdf,.eps,.ps,.png,.jpg,.jpeg}
\graphicspath{{Essai44.SWP_graphics/}{Essai44.SWP_tcache/}{Essai44.SWP_gcache/}}
\DeclareGraphicsExtensions{.pdf,.eps,.ps,.png,.jpg,.jpeg}
\graphicspath{{Essai42.SWP_graphics/}{Essai42.SWP_tcache/}{Essai42.SWP_gcache/}}
\DeclareGraphicsExtensions{.pdf,.eps,.ps,.png,.jpg,.jpeg}
\graphicspath{{Essai37.SWP_graphics/}{Essai37.SWP_tcache/}{Essai37.SWP_gcache/}}
\DeclareGraphicsExtensions{.pdf,.eps,.ps,.png,.jpg,.jpeg}
\graphicspath{{Essai35.SWP_graphics/}{Essai35.SWP_tcache/}{Essai35.SWP_gcache/}}
\DeclareGraphicsExtensions{.pdf,.eps,.ps,.png,.jpg,.jpeg}
\graphicspath{{Essai34.SWP_graphics/}{Essai34.SWP_tcache/}{Essai34.SWP_gcache/}}
\DeclareGraphicsExtensions{.pdf,.eps,.ps,.png,.jpg,.jpeg}
\graphicspath{{Essai31.SWP_graphics/}{Essai31.SWP_tcache/}{Essai31.SWP_gcache/}}
\DeclareGraphicsExtensions{.pdf,.eps,.ps,.png,.jpg,.jpeg}
\graphicspath{{Essai30.SWP_graphics/}{Essai30.SWP_tcache/}{Essai30.SWP_gcache/}}
\DeclareGraphicsExtensions{.pdf,.eps,.ps,.png,.jpg,.jpeg}
\graphicspath{{Essai28.SWP_graphics/}{Essai28.SWP_tcache/}{Essai28.SWP_gcache/}}
\DeclareGraphicsExtensions{.pdf,.eps,.ps,.png,.jpg,.jpeg}
\graphicspath{{Essai26.SWP_graphics/}{Essai26.SWP_tcache/}{Essai26.SWP_gcache/}}
\DeclareGraphicsExtensions{.pdf,.eps,.ps,.png,.jpg,.jpeg}
\graphicspath{{Essai22.SWP_graphics/}{Essai22.SWP_tcache/}{Essai22.SWP_gcache/}}
\DeclareGraphicsExtensions{.pdf,.eps,.ps,.png,.jpg,.jpeg}
\numberwithin{equation}{section}
\theoremstyle{plain}
\newtheorem{theorem}[equation]{Theorem}

\newtheorem{proposition}[equation]{Proposition}
\newtheorem{lemma}[equation]{Lemma}

\theoremstyle{definition}

\newtheorem{example}[equation]{Example}

\theoremstyle{remark}
\newtheorem{remark}[equation]{Remark}

\begin{document}
\title{Real $K$-theories }
\author{Max Karoubi  \\
Universit{\'e} Paris 7-Denis
Diderot}
\maketitle
\begin{abstract}The purpose of this short paper is to investigate relations between various real $K$-theories. In particular, we show how a real projective bundle theorem implies an unexpected relation between Atiyah's $K R$-theory and the usual equivariant $K$-theory of real vector bundles. This relation has been used recently in a new computation of the Witt group of real curves [10],
Section 4. We also interpret Atiyah's theory as a special case of twisted $K$-theory. 
\end{abstract}

\section{Short overview of classical topological $K$-theory}
Since the introduction of algebraic $K$-theory by Grothendieck [5], various versions of topological $K$-theories have emerged, essentially due to Atiyah. Historically, as a prototype of ``generalized cohomology theory'', Atiyah
and Hirzebruch [2] introduced the $K$-theory of topological complex vector bundles on a compact space $X ,$ traditionnally written $K U X) .$ There is a real analog $K O X) .\;$ When a compact group $G$ is acting continuously on $X ,$ Segal [12] defined the equivariant versions $K U_{G} (X)$ and $K O_{G} (X)$ by considering vector bundles on $X$ with a group action compatible with the action on $X ,$ We simply write $K_{G} (X)$ in a statement involving either $K U_{G} (X)$ or $K O_{G} X) .\;$ 

A more
subtle theory was introduced by Atiyah [4]. One considers a space $X$ with involution and complex vector bundles on $X ,$ with an involutive antilinear action compatible with the involution on $X ,$ Atiyah denoted this theory
by $K R (X)$ and, among other things, he showed the interest of this theory in real operator theory and in real Algebraic Geometry. This last point of
view was illustrated in recent publications [10] [11]. 

Finally, there is a less known theory starting from vector bundles which
are modules over an algebra bundle, for instance the Clifford bundle associated to a real vector bundle with a non degenerate quadratic form. This is the starting idea for the definition of the so called ``twisted $K$-theory'' with many recent publications in mathematics and physics: see e.g. [6], [13]. 

As shown some time
ago in [7], the general framework for all these theories is the concept of Banach category. For such categories $C$, one defines not only the usual Grothendieck group $K (C)$ but also ``derived functors'' $K_{i}^{t o p} C) .\;$In the same framework, for
any additive functor $\varphi  :C \rightarrow C^{ \prime }$ with an adequate continuity condition, one also define ``relative groups'' $K_{i}^{t o p} (\varphi )$ inserted in an exact sequence 

$\text{\quad \quad }K_{i +1}^{t o p} (C) \rightarrow K_{i +1}^{t o p} (C^{ \prime })$$ \rightarrow K_{i}^{t o p} (\varphi ) \rightarrow K_{i}^{t o p} (C) \rightarrow K_{i}^{t o p} (C^{ \prime })$ 

For instance, when $C$ is the category of vector bundles on $X$ and $C^{ \prime }$ the category of vector bundles on a closed subspace $Y ,$ the groups $K_{i}^{t o p} (\varphi )$ are Atiyah-Hirzebruch relative groups\protect\footnote{
We use indifferently the notation $K_{i}$ or $K^{ -i}$}$K^{ -i} (X ,Y)$. 

The interest of these relative groups will appear in the next sections. 

\section{ The real projective bundle theorems.}
Let $G$ be a compact group acting continuously on a compact space $X$ and let $V$ be an equivariant real vector bundle on $X .$ Let $W$ be an equivariant subbundle of $V ,$ $P (V)$ the real projective bundle on $X$ associated to $V$ and $P (W)$ the projective subbundle of $P( V) .\;$If $T$ is any real vector bundle, we denote by $C^{ +} (T)$ the Clifford bundle associated to $T$ with a positive quadratic form (i.e. a metric). As another definition, we denote by $\mathcal{E}^{C^{ +} (T)} (X)$ the category of vector bundles on $X$ with an action of $C^{ +}(T) .$ Finally, if a compact group
$G$ is acting on everything, we denote by $\mathcal{E}_{G}^{C^{ +} (T)} (X)$ the category of vector bundles with an intertwining action of $C^{ +} (T)$ and $G .$

\begin{theorem}
The relative group $K_{G}^{ -i} (P (V) ,P (W))$ is isomorphic to the group $K_{i -1}$ of the forgetful functor

\begin{center}
$\varphi  :\mathcal{E}_{G}^{C^{ +} (V +1)} (X) \rightarrow \mathcal{E}_{G}^{C^{ +} (W +1)} (X)$
\end{center}\par
Here the notation $T +r$ means in general the vector bundle $T$ with a trivial bundle of rank $r$ added. 
\end{theorem}

\begin{remark}
This theorem is not completely new. For $G$ trivial it has been proved in [7], Corollaire 3.2.2. However, the proof given there cannot be extended to the equivariant case,
which is important for the applications we have in mind (see e.g. Example 2.4 below). 
\end{remark}

Before giving examples, let us state
a second theorem where the notation $C^{ -} (T)$ now means the Clifford bundle of a vector bundle $T$ with a \textit{negative} non degenerate quadratic form (i.e. the opposite of a metric).

\begin{theorem}
The relative group $K_{G}^{ -i} (P (V +1) ,P (W +1))$ is isomorphic to the group $K_{i -1}$ of the forgetful functor

\begin{center}
$\Psi  :\mathcal{E}_{G}^{C^{ -} (V)} (X) \rightarrow \mathcal{E}_{G}^{C^{ -} (W)} (X)$
\end{center}\par
\begin{example}
Let us consider the usual real $K$-theory, so that $K_{G} =KO_{G} .$ Let $G =\mathbf{Z}/2$ acts on $X$ and on the trivial bundle $V =X \times \mathbf{R}$ by antipode. It is easy to see that the category $\mathcal{E}_{G}^{C^{ -} (V)} (X)$ is equivalent to the category of Real bundles on $X$ in the sense or Atiyah [4]. Therefore, the $K$-group of this category is Atiyah's group $K R( X) .$ On the other hand, if we
choose $W =0 ,$ we have $P (V +1) -P (W +1) =V =X \times \mathbf{R} ,$ so that the relative group
$K_{G} (P (V +1) ,P (W +1))$ is $K O_{G} (X \times \mathbf{R}) .\;$Therefore, Theorem 2.3 implies
the following exact sequence, apparently unknown:
\begin{flushleft}$\text{\quad \quad } \rightarrow K R (X) \rightarrow K O_{G} (X) \rightarrow K O_{G} (X \times \mathbf{R}) \rightarrow K R_{ -1} (X) \rightarrow $ (E)\end{flushleft}\par

\begin{flushleft}A similar argument with the iterated suspension of $X$ implies the more general exact sequence\end{flushleft}\par
\begin{flushleft}$\text{\quad \quad } \rightarrow K R^{ -i} (X) \rightarrow K O_{G}^{ -i} (X) \rightarrow K O_{G}^{ -i} (X \times \mathbf{R}) \rightarrow K R^{ -i +1} (X) \rightarrow $\end{flushleft}\par
\begin{remark}
\begin{flushleft}A quite different (and more involved) proof of this exact sequence is given in [10], Appendix C. \end{flushleft}\par
\end{remark}

\begin{remark}

A more concrete description of the group $K O_{G} (X \times \mathbf{R})$ is to remark that it is isomorphic to the relative group $K O_{G}( X \times D^{1} ,X \times S^{0}$). Therefore, it sits in the
middle of an exact sequence 

$\text{\quad }\qquad  \rightarrow K O_{G}^{ -1} (X) \rightarrow K O^{ -1} (X) \rightarrow K O_{G} (X \times \mathbf{R}) \rightarrow K O_{G} (X) \rightarrow $ 

For instance, if $X =Y \times S^{0} ,$ with a trivial action of
$G$ on $Y$ and the free action on $S^{0} \simeq \mathbf{Z}/2 ,$ we have $K O_{G} (X \times \mathbf{R}) \simeq K O^{ -1} Y) ,K R (X) \simeq K U (Y)$ and the exact sequence (E) is a reformulation of Bott's exact sequence for the space $Y$ (see e.g. [9] III.5.18): 

$\text{\quad \quad } \rightarrow K U (Y) \rightarrow K O (Y) \rightarrow K O^{ -1} (Y) \rightarrow K U^{1} (Y) \rightarrow $ 
\end{remark}

\begin{example}

Let $X$ be a point$ ,G$ the trivial group, $V =\mathbf{R}^{n} ,W =0.$  In this situation, the group $K_{ -1} (\Psi )$ of Theorem 2.3 is inserted in an exact sequence involving $K$-groups of Clifford algebras:
\begin{flushleft}$\text{\quad \quad }K (C^{n ,0}) \rightarrow K (C^{0 ,0}) \rightarrow K_{ -1} (\Psi ) \rightarrow K_{ -1} (C^{n ,0}) =0$\end{flushleft}\par
\end{example}

In this exact sequence, $C^{p ,q}$ denotes in general the Clifford algebra of $\mathbf{R}^{p +q}$ with the quadratic form $ -(x_{1})^{2} -\ldots  -(x_{p})^{2} +(x_{p +1})^{2} +\ldots  +(x_{p +q})^{2}$. Therefore, Theorem 2.3 implies
that the reduced $K$-theory of $R P^{n}$ is isomorphic to the cokernel of the map $K (C^{n ,0}) \rightarrow K C^{0 ,0 }) ,$ a result due to Adams [1].
Note that this result holds in real or complex $K$-theory, the Clifford algebra being understood over the relevant field.

\begin{example}
Let us consider now the relative group $K_{G} (P (V +1) ,P (V))$ which is well known to be the $K$-group $K_{G} (V)$ of the Thom space of $V .$ According to Theorem 2.1,
it is isomorphic to the $K^{1}$-group of the forgetful functor 

$\text{\quad \quad }\mathcal{E}_{G}^{C^{ +} (V +2)} (X) \rightarrow \mathcal{E}_{G}^{C^{ +} (V +1)} (X)$ 

Using Bott periodicity and the periodicity of Clifford algebras, it is not difficult to show that for $r \geq 0 ,$ this is also
the $K^{r}$-group of the forgetful functor: 

$\text{\quad \quad }\mathcal{E}_{G}^{C^{ +} (V +r +1)} (X) \rightarrow \mathcal{E}_{G}^{C^{ +} (V +r)} (X)$ 

which is a formulation of Thom's isomorphism theorem in equivariant $K$-theory [8]. However, we shall not detail this example any further since Thom's isomorphism is needed to prove Theorem 2.1, as
we shall see in next Section. 

\textit{Proof of Theorem 2,3 assuming Theorem 2.1} 

This part is purely
algebraic. According to Theorem 2.1 which we assume, the group $K_{G} (P (V +1) ,P (W +1))$ is isomorphic to the $K$-group of the forgetful functor 

$\text{\quad \quad }\mathcal{E}_{G}^{C^{ +} (V +2)} (X) \rightarrow \mathcal{E}_{G}^{C^{ +} (W +2)} (X)$ 

Note that $C^{ +} (2) =M_{2} (k) ,$ where the basic field $k =\mathbf{R}$ or $\mathbf{C}$ and `2' being an abbreviated notation\protect\footnote{
and also for the trivial bundle of rank $2.$
} for $k^{2} .$ Therefore, for any real vector bundle $T ,$ the Clifford (bundle) algebra $C^{ +} (T +2)$ is isomorphic to a $2 \times 2$ matrix algebra over $C^{ -} (T) ,\;$that is

\begin{center}
$C^{ +} (T +2) \simeq C^{ -} (T) \otimes C^{ +} (2) ,$
\end{center}\par
Indeed, we can define a map $T +2 \rightarrow C^{ -} (T) \otimes C^{ +} (2)$ by the formula $(t ,v) \mapsto t \otimes e_{1} e_{2} +1 \otimes v ,\;$$e_{1}$ and $e_{2}$ being the basis vectors of  $k^{2} .$ By the universal property of the Clifford algebra, this map induces an homomorphism $C^{ +} (T +2) \rightarrow C^{ -} (T) \otimes C^{ +} (2)$ which is an isomorphism for dimension reasons. Morita's theorem implies then that the categories $\mathcal{E}_{G}^{C^{ +} (T +2)} (T)$ and $\mathcal{E}_{G}^{C^{ -} (V)} (T)$ are naturally equivalent. As a consequence, the group $K_{G} (P (V +1) ,P (W +1))$ is isomorphic to the $K^{1}$-group of the forgetful functor $\mathcal{E}_{G}^{C^{ -} (V)} (X) \rightarrow \mathcal{E}_{G}^{C^{ -} (W)} (X) ,\;$which is essentially the formulation
of Theorem 2.3. 
\end{example}

\end{example}

\end{theorem}

\section{ Proof of Theorem 2.1}
The main ingredient to prove Theorem 2.1 is Thom's isomorphism theorem in equivariant $K$-theory. The complex analog is due to Atiyah and states that for a complex $G$-vector bundle $V$ on a $G$-space $X ,\;$the group $K U_{G} (V)$ is isomorphic to $K U_{G} X) .\;$ The general version is more
subtle and is proved in [8]. If $V$ is a real $G$-vector bundle on a $G$-space $X ,$  the group $K_{G} (V)$ is isomorphic to the Grothendieck group of the forgetful functor $\Psi  :\mathcal{E}_{G}^{C^{ +} (V +1)}(X \rightarrow \mathcal{E}_{G}^{C^{ +} (V)} (X) .\text{}\;$In particular; if $B (V)$ (resp. $S (V)$) denotes the ball bundle (resp. the sphere bundle) of $V ,$ we have an exact sequence

$ \rightarrow K_{G}(B (V) ,S (V)) \rightarrow K_{G} (B (V)) \rightarrow K_{G} (S (V)) \rightarrow K_{G}^{1}(B(V) ,S(V)) \rightarrow $ 

where $K_{G}(B (V) ,S (V)$ is identified with the $K$-group of $\Psi $. In order to prove Theorem 2.1 for the simpler case $W =0 ,\;$we have to replace the group
$G$ by $G \times \mathbf{Z}/2 ,\;$ the summand $\mathbf{Z}/2$ acting on $V$ by antipode. In that case, we may identify the group $K_{G \times \mathbf{Z}/2} (S (V))$ with the group $K_{G} (P (V))$ and the group $K_{G \times \mathbf{Z}/2} (B (V)) =K_{G \times \mathbf{Z}/2} (X)$ with the $K$-group of the category $\mathcal{E}_{G}^{C^{ +} (1)} (X)$. 

\qquad Let us look more closely at the group $K_{G \times \mathbf{Z}/2} (B (V) ,S (V))$ which is the $K$-group of the forgetful functor 

$\text{\quad \quad \quad \quad }\mathcal{E}_{G \times \mathbf{Z}/2}^{C^{ +} (V +1)} (X) \rightarrow \mathcal{E}_{G \times \mathbf{Z}/2}^{C^{ +} (V)} X) =\mathcal{E}_{G}^{C^{ +} (V +1)}$ 

Let $\eta $ represents the antipode acton of $\mathbf{Z}/2$ on $V ;$ we use $\eta $ to ``untwist'' the action of $G \times \mathbf{Z}/2$ on $C^{ +} (V +1) .\;$More precisely, by multiplying
$\eta $ with the top generator of $C^{ +} (V +1) ,\;$we show that the category
$\mathcal{E}_{G \times \mathbf{Z}/2}^{C^{ +} (V +1)} (X)$ is isomorphic to the product category $\mathcal{E}_{G \times \mathbf{Z}/2}^{C^{ +} (V)} (X) \times \mathcal{E}_{G \times \mathbf{Z}/2}^{C^{ +} (V)} (X) ,$ that is $\mathcal{E}_{G}^{C^{ +} (V +1)} (X)$ $ \times \mathcal{E}_{G}^{C^{ +} (V +1)} (X) .\text{}\;$Thanks to this isomorphism,
the previous functor may be identified with the one defined by the direct sum: 

$\text{\quad \quad }\mathcal{E}_{G}^{C^{ +} (V +1)} (X) \times \mathcal{E}_{G}^{C^{ +} (V +1)} (X) \rightarrow \mathcal{E}_{G}^{C^{ +} (V +1)} (X)$ 

We now consider the commutative diagram (up to isomorphism) 

$\begin{array}{ccc}\mathcal{E}_{G}^{C^{ +} (V +1)} (X) \times \mathcal{E}_{G}^{C^{ +} (V +1)} (X) &  \rightarrow  & \mathcal{E}_{G}^{C^{ +} (1)} (X) =\mathcal{E}_{G \times \mathbf{Z}/2} (X) \\
\downarrow  & \text{\thinspace \thinspace } & \downarrow  \\
\mathcal{E}_{G}^{C^{ +} (V +1)} (X) &  \rightarrow  &  \ast \end{array}$ 

The horizontal arrows are the algebraic analogs of the map 

$\text{\quad \quad }K_{G \times \mathbf{Z}/2} (B (V) ,S (V)) \rightarrow K_{G \times \mathbf{Z}/2} (B (V)) =K_{G \times \mathbf{Z}/2} (X)$ 

Since the $K$-theory of the first vertical functor is the $K$-theory of the category $\mathcal{E}_{G}^{C^{ +} (V +1)} (X) ,\text{}\;$by taking ``homotopy fibers''
we see that the group $K_{G \times \mathbf{Z}/2}^{ -1} (S (V)) =K_{G}^{ -1} (P (V))$ is isomorphic to the Grothendieck group of the functor 

$\text{\quad \quad }\mathcal{E}_{G}^{C^{ +} (V +1)} (X) \rightarrow \mathcal{E}_{G}^{C^{ +} (1)} (X)$ 

By considering iterated suspensions of $X$ and using Bott periodicity (and periodicity of the Clifford algebra), we deduce the theorem for $W =0.$ The general case follows from the category diagram 

$\begin{array}{ccc}\mathcal{E}_{G}^{C^{ +} (V +1)} (X) &  \rightarrow  & \mathcal{E}_{G}^{C^{ +} (1)} (X) \\
\downarrow  & \text{\thinspace \thinspace } & \downarrow  \\
\mathcal{E}_{G}^{C^{ +} (W +1)} (X) &  \rightarrow  & \mathcal{E}_{G}^{C^{ +} (1)} (X)\end{array}$ 

\section{ Relation with twisted $K$-theory}
Let us turn our attention to Atiyah's $K R$-theory of a ``real'' space $X ,\;$i.e. a space with involution.
We denote by $G$ the group \textbf{}$\mathbf{Z}/2$ and assume that $G$ is acting freely on $X$ with quotient space $Y =X/G .\;$In the case when $X =Y \times G ,$ Atiyah showed that $K R (X)$ is the usual $K U$-theory of the space $Y$ [4]. In general, let us call $L$ the real line bundle on $Y$ associated to the covering $X \rightarrow X/G$ and let us consider the Clifford bundle $C^{ -} (L)$ over $Y .\;$Note that $C^{ -} (L) \simeq 1 +L$ is the trivial complex line bundle $Y \times \mathbf{C}$ when the previous covering is trivial.

\begin{theorem}
With the above notations, the group $K R (X)$ is naturally isomorphic to the twisted $K$-theory of $Y$ associated to the Clifford algebra bundle $C^{ -} (L) .\text{}\;$

\begin{proof}
Let $E$ be a Real bundle on $X$ in the sense of Atiyah. Viewed as a real bundle, it is the pull-back of a real bundle $F$ on $Y$ so that $F =E/G .$ The complex structure
on $E$ does not come from a complex structure on $F$ (since $E$ is not a complex $G$-bundle) but it induces on the quotient by $G$ a pairing

$\text{\quad \quad }L \times F \rightarrow F$ 

and we can choose a metric on $F$ such that $F$ becomes a $C^{ -} (L)$-module. This map $E \mapsto F$ defines the required correspondance. In order to show it induces a $K$-isomorphism, we may use for instance a Mayer-Vietoris argument since the statement is true when $X =Y \times G .$ 
\end{proof}

\end{theorem}

It is interesting to describe the exact sequence in Example 2.4 

$\text{\quad \quad } \rightarrow K R (X) \rightarrow K O_{G} (X) \rightarrow K O_{G} (X \times \mathbf{R}) \rightarrow $ 

with more familiar terms.The key lemma for this is the following:

\begin{lemma}
Let $L$ be a real line bundle on a space $Y$ provided with a metric. Then the $K$-theory of the forgetful functor 

$\text{\quad \quad }\mathcal{E}^{C^{ -} (L)} (Y) \rightarrow \mathcal{E} (Y)$ 

is isomorphic to the $K$-theory of the forgetful functor associated to ``positive'' Clifford algebras 

$\text{\quad \quad }\mathcal{E}^{C^{ +} (L +2)} (Y) \rightarrow \mathcal{E}^{C^{ +} (L +1)} (Y)$

\begin{proof}
We first notice the $K$-theory equivalence between the sources and the targets. Indeed, if $E$ is vector bundle on $Y ,(L +1) \otimes E$ is a $C^{ +} (L +1)$-module, taking into account that $L \otimes L =1$ (once a metric is chosen). This correspondance induces an isomorphism between $K (Y)$ and the twisted $K$-theory associated to the algebra bundle $C^{ +} (L +1) ,\;$using Morita equivalence
and a Mayer-Vietoris argument. In the same way we proved above that the categories $\mathcal{E}^{C^{ -} (L)} (Y)$ and $\mathcal{E}^{C^{ +} (L +2)} (Y)$ are equivalent. More precisely, if $F$ is a $C^{ -} (L)$-module, we associate to it the module $F +F$ over the algebra bundle $C^{ +} (L +2) .$ The issue is now to prove that the category diagram 

$\begin{array}{ccc}\mathcal{E}^{C^{ -} (L)} (Y) &  \rightarrow  & \mathcal{E} (Y) \\
\downarrow  & \text{\thinspace \thinspace } & \downarrow  \\
\mathcal{E}^{C^{ +} (L +2)} (Y) &  \rightarrow  & \mathcal{E}^{C^{ +} (L +1)} (Y)\end{array}$ 

induces a commutative diagram between the associated $K$-groups, the vertical arrows being isomorphisms. This is checked by diagram chasing, the key remark being that if $F$ is a $C^{ -} (L)$-module, the canonical map $L \otimes F \rightarrow F$ is an isomorphism.

\begin{proposition}
Let $X$ be a free $G$-space (with $G =\mathbf{Z}/2)$ and $Y =X/G .\;$ Then the general exact sequence

$ \rightarrow K O_{G}^{ -1} (X \times \mathbf{R}) \rightarrow K R (X) \rightarrow K O_{G} (X) \rightarrow K O_{G} (X \times \mathbf{R}) \rightarrow $ 

is isomorphic to the exact sequence associated to twisted $K$-groups 

$ \rightarrow K (\Psi ) \rightarrow K (\mathcal{E}^{C^{ -} (L)} (Y)) \rightarrow K (\mathcal{E} (Y)) \rightarrow K^{1} (\Psi ) \rightarrow $ 

where $\Psi $ is the forgetful functor $\mathcal{E}^{C^{ -} (L)} (Y) \rightarrow \mathcal{E} (Y) .$

\begin{proof}
Let $L$ be the real line bundle associated to the covering $X \rightarrow Y .\;$According to the lemma, the
$K$-theory of $\Psi $ is isomorphic to the $K$-theory of the functor $\Theta  :\mathcal{E}^{C^{ +} (L +2)} (Y) \rightarrow \mathcal{E}^{C^{ +} (L +1)} (Y) .\;$On the other hand, the $K$-theory Thom's isomorphism proved in [7] p. 210 shows that $K (\Theta )$ is isomorphic to $KO (L +1) =KO^{ -1}(L) =K O_{G}^{ -1} (X \times \mathbf{R}) .\;$ 
\end{proof}

\end{proposition}

\end{proof}

\end{lemma}

\begin{center}
\textbf{REFERENCES}
\end{center}\par
[1] J.-F. Adams. Vector fields on the sphere. Ann. Math. 75 (1962),
603-632. 

[2] M.-F. Atiyah and F. Hirzebruch. Vector bundles and homogeneous spaces. Symposia in Pure Math. Amer. Math. Society 3 (1961),
7-38. 

[3] M.F. Atiyah, R. Bott and A. Shapiro. Clifford modules. Topology 3 (1964), 3-38 

[4] M.F. Atiyah. $K$-theory and reality. Quarterly Math. J. Oxford (2) (1966), 367-386. 

[5] A. Borel et J.-P. Serre. Le th{\'e}or{\`e}me
de Riemann-Roch (d'apr{\`e}s Grothendieck). Bull. Soc. Math. France 86 (1958), 97-136. 

[6] P. Donovan and M. Karoubi. Graded Brauer
groups and $K$-theory with local coefficients. Inst. Hautes Etudes Sci. Pub. Math. N
38 (1970), 5-25. 

[7] M. Karoubi. Alg{\`e}bres de Clifford et $K$-th{\'e}orie. Ann. Sci. Ec. Norm Sup. (4) (1968), 161-270. 

[8] M. Karoubi. Sur le theor{\`e}me de Thom
en $K$-th{\'e}orie {\'e}quivariante. Springer Lecture Notes N 136 (1970),
187-253. 

[9] M. Karoubi. K-theory. An introduction.Grundlheren der Math. N
226 (1978) 

[10] M. Karoubi, M. Schlichting and C. Weibel. The Witt group of real algebraic varieties. J. Topol. 9 (2016), 1257-1302.

[11] M. Karoubi and C. Weibel. The Witt group of real surfaces. , Contemporary Math. 749, Amer. Math. Society (2020), 157-193. 

[12] G. Segal. Equivariant $K$-theory. Inst. Hautes Etudes Sci. Pub. Math. N 34 (1968), 129-151. 

[13] E. Witten. $D$-branes and $K$-theory. J. High Energy Phys. 12, Nr. 19 (1998). 
\end{document}